\documentclass[12pt]{article}
\usepackage[margin=1in]{geometry}
\usepackage{graphicx} 

\usepackage{cite}
\usepackage{amssymb}
\usepackage{amsthm}
\usepackage{amsmath}

\usepackage{bigdelim}

\usepackage{comment}
\usepackage{xcolor}
\usepackage{hyperref}

\usepackage{nicematrix}

\theoremstyle{plain}
\newtheorem{thm}{Theorem}[section]

\newtheorem{prop}[thm]{Proposition}
\newtheorem{cor}[thm]{Corollary}

\newtheorem*{thm*}{Theorem}
\newtheorem*{lemma*}{Lemma}
\newtheorem*{prop*}{Lemma}
\newtheorem*{cor*}{Corollary}
\newtheorem*{conj*}{Conjecture}

\theoremstyle{definition}

\newtheorem*{defn*}{Definition}
\newtheorem{ex}[thm]{Example}

\theoremstyle{remark}
\newtheorem{rmk}[thm]{Remark}

\DeclareMathOperator{\aff}{Aff}
\renewcommand{\sl}{{\rm SL}}
\newcommand{\gl}{{\rm GL}}

\usepackage{tikz}
\tikzstyle{vertex}=[circle, draw, inner sep=0pt, minimum size=4pt, fill=black]
\newcommand{\vertex}{\node[vertex]}

\title{The pure condition for incidence geometries}
\author{Daniel Irving Bernstein and Signe Lundqvist}
\date{}

\begin{document}

\maketitle

\begin{abstract}

The space of \emph{parallel redrawings} of an incidence geometry $(P,H,I)$ with an assigned set of normals is the set of points and hyperplanes in $\mathbb{R}^d$ satisfying the incidences given by $(P,H,I)$, such that the hyperplanes have the assigned normals. 

In 1989, Whiteley characterized the incidence geometries that have $d$-dimensional realizations with generic hyperplane normals such that all points and hyperplanes are distinct. However, some incidence geometries can be realized as points and hyperplanes in $d$-dimensional space, with the points and hyperplanes distinct, but only for specific choices of normals. Such incidence geometries are the topic of this article.

In this article, we introduce a pure condition for parallel redrawings of incidence geometries, analogous to the pure condition for bar-and-joint frameworks, introduced by White and Whiteley. The $d$-dimensional pure condition of an incidence geometry $(P,H,I)$ imposes a condition on the normals assigned to the hyperplanes of $(P,H,I)$
required for $d$-dimensional realizations of $(P,H,I)$ with distinct points.
We use invariant theory to show that the $d$-dimensional pure condition of $(P,H,I)$ is a bracket polynomial.
We will also explicitly compute the pure condition as a bracket polynomial for some examples in the plane.
\end{abstract}

\section{Introduction}

Given a finite sets of points $P \subseteq \mathbb{R}^d$ and a finite set $H$ consisting of affine hyperplanes in $\mathbb{R}^d$, the corresponding \emph{incidence geometry} is the triple $(P,H,I)$ where
$I \subseteq P \times H$ consists of the pairs $(p,h)$ such that $p$ lies in $H$.
As with any combinatorial abstraction of a geometric object, one can go the other direction,
and ask about the configuration of points and hyperplanes that have some particular prescribed incidence geometry.
In this paper, we are particularly interested in constraints that the combinatorics of an incidence geometry $(P,H,I)$ place on the normal vectors of hyperplanes in a realization.

Given an incidence geometry $(P,H,I)$ and a function $n: H\rightarrow \mathbb{R}^d$ assigning a normal direction to each hyperplane in $H$,
the corresponding space of $d$-dimensional \textit{parallel redrawings} is the space of realizations of $(P,H,I)$ as points and hyperplanes in $\mathbb{R}^d$ where the hyperplane $h \in H$
has normal direction $n(h)$.

It is also natural to consider \textit{polyhedral scenes}, which are the dual to parallel redrawings. The fundamental problem in scene analysis is reconstructing a $3$-dimensional polyhedral cap from its projection to $2$-dimensional space \cite{SUGIHARA, Plane_Matroid, Scene_analysis_frameworks}. More specifically, polyhedral scenes are realizations of an incidence geometry as points and hyperplanes in $\mathbb{R}^d$ such that projection of the first $d-1$ coordinates to $\mathbb{R}^{d-1}$ is specified. The results of this paper all have corresponding dual results for polyhedral scenes, that can be proven in the same way. 

There is a classical theory dating back to Maxwell, which gives a correspondence from parallel redrawings of polyhedra to stresses in a reciprocal diagram that can be obtained from the normals of the faces of the polyhedra \cite{Maxwell, C_W_Spaces_of_Stresses, C_W_Plane_self_stresses}. In the plane, there is also a one-to-one correspondence between the parallel redrawings of a bar-joint frameworks and its infinitesimal motions \cite{discrete_matroids, Plane_Matroid}. 

More recently, (a slight generalization of) parallel redrawings were used to understand the Jacobian ideal of hyperplane arrangements \cite{Hyperplanes}. 

Every space of parallel redrawings contains the \emph{trivial realizations}, which are realizations where all the points coincide. The incidence geometries whose spaces of parallel redrawings are nontrivial for generic normal vector direction functions $n: H \rightarrow \mathbb{R}^d$ are the non-spanning independent sets of a particular matroid studied in~\cite{Plane_Matroid}.
For spanning sets of the matroid, there may exist choices of normal vector direction functions $n: H\rightarrow \mathbb{R}^d$
such that the corresponding spaces of nontrivial $d$-dimensional parallel redrawings are nonempty, but such normals are not generic.
As we will show, they are the points lying in a particular algebraic hypersurface defined by a bracket polynomial.

This generalizes an analogous situation for bar and joint frameworks, the main object of study in classical rigidity theory.
The pure condition for bar-and-joint frameworks, introduced by White and Whiteley, is an algebraic condition on the coordinates of bar-joint frameworks realizing a generically minimally rigid graph that are infinitesimally flexible \cite{pure_condition:1}.
The pure condition is obtained by choosing a ``tie-down" of the framework to eliminate the Euclidean motions, and then taking the determinant of the rigidity matrix. White and Whiteley prove that the pure condition is independent of the chosen tie-down, and depends only on the underlying combinatorial graph \cite{pure_condition:1}. Furthermore, the pure condition can be expressed as a bracket polynomial \cite{pure_condition:1}. Pure conditions which describe when structures exhibit non-generic behaviour have also been described for bar-and-body frameworks and body-cad frameworks \cite{pure_condition:2, pure_condition:body-cad}. 

In this paper, we introduce a pure condition for realizations of incidence geometries. In this setting, the pure condition is an algebraic condition on the normals assigned to the hyperplanes that needs to be satisfied for the incidence geometry to have realizations with those line slopes and all points distinct. The pure condition can be obtained by pinning a point to remove translations, and then considering the determinant of the matrix which has as its kernel the parallel redrawings of the incidence geometry. We will prove that the pure condition does not depend on the pinned point.
We will use invariant theory to prove that it can be expressed as a bracket polynomial. We will also explicitly compute the pure condition as a bracket polynomial for some examples.

\section{Preliminaries}

\subsection{Notation}
Given a finite set $S$,
the $\mathbb{R}$-vector space whose coordinates are indexed by $S$ is denoted $\mathbb{R}^S$.
For a positive integer $d$,
the set of $\mathbb{R}$-matrices with $d$ rows whose columns are indexed by $S$ will be denoted $\mathbb{R}^{d\times S}$.
Given $M \in \mathbb{R}^{d\times S}$, we denote the column of $M$ corresponding to $e \in S$ by $M_e$.
The special and general linear groups of a vector space $V$ will be denoted $\sl(V)$ and $\gl(V)$.

\subsection{Parallel redrawings and associated matroids}

An \textit{incidence geometry} is a triple of sets $(P,H,I)$, where $P$ is a set of \textit{points}, $H$ is a set of \textit{hyperplanes} and $I \subseteq P \times H$ is a set of incidences between the points and hyperplanes.
One can view an incidence geometry as a hypergraph where the sets of points incident to a common hyperplane are the hyperedges.
For each integer $d \ge 1$, the set of affine hyperplanes in $\mathbb{R}^d$ will be denoted $\aff(d)$.
Given an incidence geometry $(P,H,I)$,
a \emph{$d$-dimensional realization} of $(P,H,I)$ is a pair $(\rho,\sigma)$ consisting of functions
$\rho: P \rightarrow \mathbb{R}^d$ and $\sigma: H \rightarrow \aff(d)$
such that $\rho(p)$ lies in $\sigma(h)$ whenever $(p,h) \in I$.
\begin{ex}\label{ex: 2d incidence geometry example}
    Define $P = \{p_0,\dots,p_6\}$ and $H = \{h_0,\dots,h_5\}$ and 
    \begin{align*}
        I = \{(p_0,h_0), (p_1, h_0), (p_2, h_0), \\
            (p_0,h_1), (p_3, h_1), (p_4, h_1), \\
            (p_0,h_2), (p_5, h_2), (p_6, h_2), \\
            (p_1,h_3), (p_3,h_3), (p_6,h_3), \\
            (p_2,h_4), (p_3,h_4), (p_5,h_4), \\
            (p_2,h_5), (p_4,h_5), (p_6,h_5) \}
    \end{align*}
    A two-dimensional realization of the incidence geometry $(P,H,I)$ can be found in Figure \ref{fig:Non-Fano}.
    \begin{figure}
    \[
        \begin{tikzpicture}
            \vertex (1) at (0,0)[label=left:$p_6$]{};
            \vertex (2) at (2,0)[label=right:$p_2$]{};
            \vertex (3) at (1,2)[label=above:$p_0$]{};
            \vertex (4) at (1,0)[label=below:$p_4$]{};
            \vertex (5) at (1/2,1)[label=left:$p_5$]{};
            \vertex (6) at (3/2,1)[label=right:$p_1$]{};
            \vertex (7) at (1,43/64)[label=right:$p_3$]{};
            \path[thick]
                (1) edge (3)
                (1) edge (2)
                (2) edge (3)
                (1) edge (6)
                (2) edge (5)
                (3) edge (4)
            ;
            \node at (0,1/2){$h_2$};
            \node at (1/2,-1/4){$h_5$};
            \node at (5/8,3/8){$h_3$};
            \node at (3/4,1){$h_4$};
            \node at (1,6/4){$h_1$};
            \node at (2,1/2){$h_0$};
        \end{tikzpicture}
    \]
    \caption{A two-dimensional realization of an incidence geometry.}
    \label{fig:Non-Fano}
\end{figure}
\end{ex}

Given a function $n: H \rightarrow \mathbb{R}^d$ associating a normal vector to each element of $H$,
a $d$-dimensional \emph{parallel redrawing} of $(P,H,I;n)$ is a $d$-dimensional realization $(\rho,\sigma)$ such that $n(h)$ is a normal vector of $\sigma(h)$ for each $h \in H$.
In other words, a parallel redrawing of $(P,H,I)$ with a prescribed set of normals consists of functions $\iota: H \rightarrow \mathbb{R}$ and $x: P \rightarrow \mathbb{R}^d$ such that
\begin{equation}
    n(h)\cdot x(p)+\iota(h)=0 \qquad {\rm for \ all } \qquad (p,h) \in I.
    \label{parallel}
\end{equation}

\begin{rmk}
    A common assumption when working with parallel redrawings is that the last coordinate of $n(h)$ is $1$ for all hyperplanes, which says that the hyperplane is not parallel to the hyperplane $x_d=0$ in $\mathbb{R}^d$. If we assume that the normals are generic, i.e. that the coordinates of $n(h)$ are algebraically independent over $\mathbb{Q}$, then in particular the last coordinate of $n(h)$ has to be non-zero, so we might assume that it is $1$, by multiplying equation \ref{parallel} by a constant. 
    
    In later sections, we will not make the assumption that the last coordinate of $n(h)$ is $1$.
    If an incidence geometry has proper realizations where one of the hyperplanes is parallel to the hyperplane $x_d=0$,
    then one can rotate to obtain a proper realization such that none of the hyperplanes are parallel to the hyperplane $x_d=0$.
\end{rmk}

Finding the parallel redrawings of an incidence geometry $(P,H,I)$ with an associated normal vector $n(h)$ to each hyperplane amounts to solving a system of equations of the form~\eqref{parallel},
one for each incidence in $I$. 
This system has $|H|+d|P|$ indeterminates, one for each $h \in H$ from $\iota(h)$, and $d$ for each $p \in P$ from $x(p)$. 

Let the matrix $S \in \mathbb{R}^{d\times H}$ be the matrix with columns given by the normal vectors of the hyperplanes. Let $M_S(P,H,I)$ denote the $|I| \times (|H|+d|P|)$ coefficient matrix of the system of equations of the form~\eqref{parallel}.
Each $p \in P$ indexes $d$ columns and each $h \in H$ indexes a single column.
The row of $M_S(P,H,I)$ corresponding to the incidence $(p, h)$ is

  \setcounter{MaxMatrixCols}{20}
  \[\begin{bNiceMatrix}[first-row,first-col]
         &&&&h&&&&&p&&&&&\\
		  &0&...&0&1&0&...&0&n(h)_1&\dots&n(h)_d&&0&...&0
        \end{bNiceMatrix}\]
 where the entries of $n(h)$ are placed at the columns indexed by $p$, and the $1$ is placed at the column indexed by $h$.

 There are realizations of any incidence geometry $(P,H,I)$ with any choice of normals $n$ given by assigning the same point coordinates to all elements of $P$. Such realizations are called \textit{trivial}.
 A realization is \emph{proper} if no two points nor two hyperplanes coincide.
 Not all incidence geometries have proper realizations. The Fano plane is an example of an incidence geometry that does not have proper realizations for \textit{any} choice of hyperplane normals.

 \begin{rmk}
     Given an incidence geometry $(P,H,I)$ such that any pair of points are incident to at most one common hyperplane, there is a corresponding rank $3$ matroid $\mathcal{M}(P,H,I)$ with ground set $P$ and bases given by the triples of points that are not collinear. In particular, $\mathbb{R}$-representations are proper realization of $(P,H,I)$ in $\mathbb{R}^2$.
 \end{rmk}

There is always a $d$-dimensional space of parallel redrawings of any realization of an incidence geometry in $d$-dimensional space,
including trivial realizations, given by translations.
 
 Let $S=(P,H,I)$ be an incidence geometry. The \textit{$d$-plane matroid} is the matroid defined on the set of incidences $I$ where a set $I' \subseteq I$ is said to be independent if $|I''| \leq |H(I'')| + d|P(I'')|-d$ for all non-empty subsets $I'' \subseteq I'$,
 where $P(I'')$ and $H(I'')$ denote the projections of $I'' \subseteq P \times H$ onto $P$ and $H$.
 That this is indeed a matroid follows from the next theorem.

 \begin{thm}[Whiteley \cite{Plane_Matroid}]\label{thm: whiteley plane matroid}
     Let $(P,H,I)$ be an incidence geometry.
     The rows of $M_S(P,H,I)$ are independent for any choice of generic normals $n$ if and only if $I$ is independent in the $d$-plane matroid.
     \label{picture_theorem}
 \end{thm}

A corollary to Theorem \ref{picture_theorem} is that an incidence geometry $(P,H,I)$ has proper parallel redrawings in dimension $d$ for any generic choice of normals if and only if $|I''| \leq |H(I'')|+d|P(I'')|-(d+1)$ for all subsets $I'' \subseteq I$ with $|I''| \geq 2$.
Note that this condition excludes sets of incidences that span the $d$-plane matroid.
Some incidence geometries that do not satisfy this condition will have proper parallel redrawings for some non-generic choices of normals -- consider Example \ref{ex: 2d incidence geometry example}.
It has $7$ points, $6$ hyperplanes and $18$ incidences, so it satisfies $|I|=|H|+2|P|-2$.
In fact, it is a basis in the $2$-plane matroid, which means that it has only trivial realizations in the plane for generic normals.
However, as can be seen in Example \ref{ex: 2d incidence geometry example}, it does have proper realizations for some choice of normals.
In this paper, we will develop an approach for finding such choices of normals for incidence geometries that span the $d$-plane matroid.
 
\subsection{Invariant theory}
We now introduce the minimum necessary background from invariant theory.
Let $V$ be an $\mathbb{R}$-vector space and
let $\mathbb{R}[V]$ denote the ring of polynomial functions on $V$,
and let $\Gamma$ be a subgroup of $\gl(V)$.
The \emph{invariant ring of $\Gamma$} is
\[
    \mathbb{R}[V]^\Gamma := \{f \in \mathbb{R}[V] : f = f \circ \pi {\rm \ for \ all \ } \pi \in \Gamma\}.
\]
When $V$ comes naturally equipped with a basis $e_1,\dots,e_k$,
we will represent $\mathbb{R}[V]$ as $\mathbb{R}[x_1,\dots,x_k]$.
In this case, given $f \in \mathbb{R}[x_1,\dots,x_k]$,
the value of $f$ at $\sum_{i = 1}^k a_i e_i$ is $f(a_1,\dots,a_k)$.

Let $d,n$ be positive integers with $n \ge d$
and let $X$ be the $d\times n$ matrix whose $ij$ entry is the indeterminate $x_{ij}$;
the corresponding real polynomial ring is denoted $\mathbb{R}[X]$.
For each tuple $(i_1,\dots,i_d)$ such that $i_j < i_k$ whenever $j<k$,
the determinant of the $d\times d$ submatrix of $X$ obtained by restricting to the columns indexed by $i_1,\dots,i_d$ is denoted $[i_1,\dots,i_d]$.
If the $i_j$s are out of order, then $[i_1,\dots, i_d]$ denotes this determinant, times the sign of the permutation needed to put the $i_j$s in order.
In the statement of the following theorem, we are viewing $\sl(\mathbb{R}^{d})$ as a subgroup of $\gl(\mathbb{R}^{d\times n})$,
acting column-wise.

\begin{thm}[First fundamental theorem of invariant theory]\label{thm: first fundamental theorem of invariant theory}
    $\mathbb{R}[X]^{\sl(\mathbb{R}^{d})} = \mathbb{R}[[i_1,\dots,i_d] : 1 \le i_1 < \dots < i_d \le n]$.
\end{thm}
\begin{proof}
    The analogous result for the field $\mathbb{C}$ can be found in~\cite[Theorem 3.2.1]{sturmfels2008algorithms}.
    The same proof technique works for arbitrary fields of characteristic zero, see also~\cite[Theorem~6.14]{neusel2007invariant}.
\end{proof}
Elements of the invariant ring in Theorem~\ref{thm: first fundamental theorem of invariant theory} are called \emph{bracket polynomials}.
Theorem~\ref{thm: first fundamental theorem of invariant theory} enables us to prove that a given polynomial $f$ is a bracket polynomial by showing that $f = f \circ \pi$ for each $\pi \in \sl(\mathbb{R}^{d \times n})$.

\section{The pure condition}
In this section, we will define a pure condition for parallel redrawings of incidence geometries, analogous to the pure condition for bar-and-joint frameworks defined by White and Whiteley \cite{pure_condition:1}.
Similarly to the case of bar and joint frameworks, we will prove that the pure condition for incidence geometries
can be described at the zero-set of a bracket polynomial.

Let $(P,H,I)$ be an incidence geometry.
For each positive integer $d$ we define a variety associated with $(P,H,I)$ as follows.
Let $X$ be a matrix of indeterminates with $d$ rows whose columns are indexed by $P$,
and let $y$ be a vector of indeterminates whose entries are indexed by $H$.

To each incidence $(p,h) \in I$ and each matrix $S \in \mathbb{R}^{d\times H}$ of normal vectors,
we associate the linear functional $A_S^{d,p,h}$ on $\mathbb{R}^{|H|+d|P|}$, defined as follows
\begin{equation}\label{eq: point hyperplane incidence}
    A_S^{d,p,h}(X,y) = \langle S_h, X_p\rangle + y_h.
\end{equation}
Then $M_S(P,H,I) \in \mathbb{R}^{I \times (|H|+d|P|)}$ is the matrix whose $(p,h)$ row is the coefficient vector of $A_S^{d,p,h}$. The kernel of $M_S(P,H,I)$ gives the set of realizations of $(P,H,I)$ such that the normal vectors of $H$ are given by $S$.

The kernel of $M_S(P,H,I)$ will always have at least $d$ dimensions since any translation of a $d$-realization of $(P,H,I)$
is another $d$-realization with the same normal vectors.
This is an extra $d$ dimensions of redundant information,
which we can remove by restricting the kernel of $M_S(P,H,I)$ to only those realizations that leave some particular $p \in P$ fixed at the origin.
Indeed, such configurations comprise the kernel of a new matrix $M_S^p(P,H,I)$
which we define as follows:
for each of the $d$ columns corresponding to $p$,
add a row that is all zeros aside from a $1$ at that column.
The dimension of the set of such realizations should not depend on the choice of point $p \in P$ that gets fixed to the origin.
Indeed, it does not as the next proposition shows.

\begin{prop}\label{rank_independent}
    Let $(P,H,I)$ be an incidence geometry, and let $S$ be some choice of normals.
    Then the rank of $M^p_S(P,H,I)$ is independent of the choice of fixed point $p$. 
\end{prop}
\begin{proof}
    Fix $p\in P$ and let $v \in \mathbb{R}^{|H|+d|P|}$ be such that $M_{S}^p(P,H,I)v = 0$.
    Given $q \in P$ and $h \in H$, the $d$ coordinates of $v$ corresponding to $q$ will be denoted $X_q$
    and the coordinate of $v$ corresponding to $h$ will be denoted $y_h$.
    As $p$ is fixed, $X_p = 0$.
    Fix some $q \in P$.
    For each $w \in P$, define $\hat X_w := X_w-X_q$ and for each $h \in H$ define $\hat y_h = y_h + \langle S, X_q \rangle$
    and let $\hat v$ be the corresponding element of $\mathbb{R}^{d|P|+|H|}$.
    Then the map $v \mapsto \hat v$ is a bijection from the kernel of $M_S^p(P,H,I)$ to the kernel of $M_S^q(P,H,I)$. It follows that the dimension of the kernel of $M^p_S(P,H,I)$ is independent of 
    the choice of fixed point $p$, so the rank must also be independent of the choice of fixed point $p$.
\end{proof}

Now assume that $(P,H,I)$ is an incidence geometry such that $I$ is a basis of the $d$-plane matroid,
i.e.~that it satisfies $|I| = |H|+d|P| - d$ and that
$|I'| \le d|P'|+|H'|-d$ for every $I' \subseteq I$ where $P',H'$ are the points and hyperplanes appearing in $I'$.
In this case $M_S^p(P,H,I)$ is a square matrix.
For generic fixed $S$, the kernel of $M_S(P,H,I)$ is $d$-dimensional, consisting only of trivial realizations that place each $p \in P$ at the same point,
and the kernel of $M_S^p(P,H,I)$ consists only of the origin.
For arbitrary $S$, the kernel of $M_S(P,H,I)$ has dimension at least $d+1$ if and only if $M^p_S(P,H,I)$ is singular.
In this case, the kernels of $M_S(P,H,I)$ and $M^p_S(P,H,I)$ both contain nontrivial $d$-realizations of $(P,H,I)$ with normal vectors $S$.
Thus for a given $S$, we can test if $(P,H,I)$ contains nontrivial $d$-realizations of $(P,H,I)$ with normal vectors in $S$
by seeing if the determinant of $M^p_S(P,H,I)$ is zero.
Indeed, by treating the entries of $S$ as indeterminates, we can view the determinant of $M_S^p(P,H,I)$ as a polynomial
It is square, so we can compute its determinant, which is a polynomial in the entries of $S$.
The vanishing locus of this polynomial gives us the normal vectors for which nontrivial $d$-realizations of $(P,H,I)$ exist.
    
For bar-joint frameworks, body-bar frameworks and body-cad frameworks, the determinant of the rigidity matrix of a tied-down framework can be factored as a factor that depends on the tie-down, and a critical factor that depends only on the underlying combinatorial structure \cite{pure_condition:1, pure_condition:body-cad}. The next proposition says that in the case of parallel redrawings, the determinant of the pinned matrix depends on the pinned point only up to a constant factor.

\begin{prop}\label{prop: pure condition doesn't depend on choice of pinned point}
    Let $(P,H,I)$ be an incidence geometry such that $I$ is independent in the $d$-plane matroid.
    Given $p,q \in P$ there exists a nonzero $\lambda \in \mathbb{R}$ such that for all $S \in \mathbb{R}^{d\times H}$
    \[
        \det M_S^p(P,H,I) = \lambda \det M_S^q(P,H,I).
    \]
\end{prop}
\begin{proof}
    For each $p \in P$, define
    \[
        f_p(S) = \det(M_S^p(P,H,I)).
    \]
    We first show that $f_p$ is multilinear and has degree that is independent of the choice of $p \in P$.
    Indeed, each monomial of $f_p$ is of the following form for some permutation $\sigma$ and some coefficient $c \in \mathbb{R}$
    \begin{equation}
        c \Pi_{i=1}^{|H|+d|P|}M_S^p(P,H,I)_{i,\sigma(i)}.
        \label{monomial}
    \end{equation}
    It suffices to show that if $c \neq 0$, then this monomial has degree $d(|P|-1)$.
    Indeed, $M_S(P,H,I)$ has $|H|$ columns of $\{0,1\}$-vectors
    and $d|P|$ columns whose entries are either $0$ or an indeterminate.
    One obtains $M_S^p(P,H,I)$ from $M_S(P,H,I)$ by adding the $d\times d$ identity submatrix
    below the columns corresponding to $p$, and $d$ zeros below every other column.
    Thus the determinant of $M_S^p(P,H,I)$ is equal to the determinant of the matrix obtained from $M_S(P,H,I)$ by deleting the $d$ columns corresponding to $p$,
    and in this matrix, each of the $d(|P|-1)$ columns corresponding to the other elements of $P$ contributes one to the degree of the determinant.

    Now let $p,q \in P$.
    Then $f_p$ and $f_q$ are multilinear of the same degree.
    Proposition~\ref{rank_independent} implies that $f_p$ and $f_q$ have the same vanishing locus.
    Therefore $f_p = \lambda f_q$ for some nonzero $\lambda \in \mathbb{R}$.
\end{proof}

Let $(P,H,I)$ be an incidence geometry that is a basis of the $d$-plane matroid, let $X$ be a $d\times |H|$ matrix of indeterminates, and let $p \in P$.
The \emph{pure condition} of $(P,H,I)$ is the following polynomial constraint on the entries of $X$
\[
    \det(M_X^p(P,H,I)) = 0.
\]
The columns of a matrix $S \in \mathbb{R}^{d\times H}$ are the normal vectors of a nontrivial realization of $(P,H,I)$
if and only if $S$ satisfies the pure condition.
Proposition~\ref{prop: pure condition doesn't depend on choice of pinned point} implies that the pure condition does not depend on choice of $p$.

Our next order of business is to show that the pure condition can be expressed as a bracket polynomial,
just like in the case of bar and joint frameworks.
We do this via invariant theory.
Given a finite set $H$,
the \emph{diagonal action of $\sl(\mathbb{R}^d)$ on $\mathbb{R}^{d\times H}$}
is given by applying the canonical action of $\sl(\mathbb{R}^d)$ on $\mathbb{R}^d$ column-wise to matrices in $\mathbb{R}^{d\times H}$.
Given $A \in \sl(\mathbb{R}^d)$ and $S \in \mathbb{R}^{d\times H}$,
the image of $S$ under the diagonal action of $A$ is simply $AS$.

\begin{prop}\label{prop: parallel redrawing polynomials are SL invariant}
    Let $d \ge 2$ be an integer and let $(P,H,I)$ be an incidence geometry satisfying $|I| = d|P|+|H|-d$.
    Let $p \in P$ and let $S \in \mathbb{R}^{d\times H}$.
    Then the determinant of $M_S^p(P,H,I)$ is invariant with respect to the diagonal action of $\sl(\mathbb{R}^d)$ on $S$,
    i.e.~if $A \in \sl(\mathbb{R}^d)$ and $S \in \mathbb{R}^{d\times |H|}$ then
    \[
        \det(M_S^p(P,H,I)) = \det(M_{A S}^p(P,H,I)).
    \]    
\end{prop}
\begin{proof}
    Fix $p \in P$ and $S \in \mathbb{R}^{d\times H}$ and $A \in \sl(\mathbb{R}^d)$.
    Let ${\rm Id}_{H}$ denote the $|H| \times |H|$ identity matrix and
    let $\hat A$ denote the $(|H|+d|P|-d) \times (|H|+d|P|-d)$ block-diagonal matrix defined as follows
    \[
        \hat A :=
        \begin{pmatrix}
            {\rm Id}_H & 0 & \dots & 0 & 0\\
            0 & A^T & \dots & 0 & 0\\
            0 & 0 & \dots & A^T & 0\\
            0 & 0 & \dots & 0 & A^T
        \end{pmatrix}.
    \]
    Then $\det(\hat A) = 1$. The result now follows since $M_S^p(P,H,I) \hat A = M_{AS}^p(P,H,I)$.
\end{proof}

\begin{cor}\label{cor: these are bracket polynomials}
    Let $d \ge 2$ be an integer and let $(P,H,I)$ be an incidence geometry satisfying $|I| = d|P|+|H|-d$, let $p \in P$,
    and let $X$ be a $d\times |H|$ matrix of indeterminates.
    Then $\det(M_X^p(P,H,I))$ is a bracket polynomial in $\mathbb{R}[X]$.
\end{cor}
\begin{proof}
    This follows from Proposition~\ref{prop: parallel redrawing polynomials are SL invariant} via Theorem~\ref{thm: first fundamental theorem of invariant theory}.
\end{proof}

\section{Two-dimensional examples}
In this section we explicitly compute the pure condition for a few examples in the case that $d = 2$.
Unless otherwise stated, all configurations $(P,H,I)$ will be assumed to satisfy the property that $I$ is a basis of the $2$-plane matroid.
The structure of the matrix $M^p_S(P,H,I)$ can be exploited to explicitly compute $\det(M^p_S(P,H,I))$ as a bracket polynomial.
We will fix orderings $P = \{p_0,\dots,p_{|P|-1}\}$ and $H = \{h_0,\dots,h_{|H|-1}\}$
and assume that the pinned point $p$ is equal to $p_{|P|-1}$.

Let $n: H \rightarrow {\rm Aff}(2)$ and let $(f_i,g_i)$ denote the normal vector of $n(h_i)$.
Suppose that $\{p_0, p_1,...,p_k\}$ is the set of points incident to $h_0$.
Then the row-submatrix of $M^p_S(P,H,I)$ corresponding to the incidences $\{(h_0,p_0), (h_0,p_1),...,(h_0, p_k)\}$ looks as follows:
\[
\begin{array}{ccccccccccccccccccccccccc}
     & & h_0 & \dots & h_{|H|-1} & p_0 & p_1 & \dots & p_k & p_{k+1} & \dots & p_{|P|-2}& & \\
   (h_0,p_0) & \ldelim({5}{0.5em} & 1 &0 &0 &f_0 \quad g_0 &0 \quad 0 &0 \quad 0 &0 \quad 0 & 0 \quad 0&\dots &0 \quad 0& \rdelim){5}{0.5em} \\
   (h_0,p_1) && 1 & 0 & 0 & 0 \quad 0& f_0 \quad g_0 &0 \quad 0 &0 \quad 0 &0 \quad 0&\dots &0 \quad 0&\\
   \vdots && \vdots & \vdots & \vdots & \vdots \quad \vdots & \vdots \quad \vdots & \vdots \quad \vdots & \vdots \quad \vdots & \vdots \quad \vdots &\dots& \vdots \quad \vdots&\\

   (h_0,p_k) && 1 & 0 & 0 & 0 \quad 0 &0 \quad 0 &0 \quad 0 & f_0 \quad g_0 & 0 \quad 0 &\dots & 0 \quad 0 &
\end{array}
\]
\bigskip

Subtracting the row corresponding to the incidence $(h_0,p_0)$ from the rows corresponding to the remaining incidences $(h_0,p_i)$ for $1 \leq i \leq k$ yields the row submatrix

\[
\begin{array}{ccccccccccccccccccccccccc}
     & & h_0 & \dots & h_{|H|-1} & p_0 & p_1 & \dots & p_k & p_{k+1} & \dots & p_{|P|-2}& & \\
   (h_0,p_0) & \ldelim({5}{0.5em} & 1 &0 &0 &f_0 \quad g_0 &0 \quad 0 &0 \quad 0 &0 \quad 0 & 0 \quad 0&\dots &0 \quad 0& \rdelim){5}{0.5em} \\
   (h_0,p_1) && 0 & 0 & 0 & -f_0 \quad -g_0 & f_0 \quad g_0 &0 \quad 0 &0 \quad 0 &0 \quad 0&\dots &0 \quad 0&\\
   \vdots && \vdots & \vdots & \vdots & \vdots \quad \vdots & \vdots \quad \vdots & \vdots \quad \vdots & \vdots \quad \vdots & \vdots \quad \vdots &\dots& \vdots \quad \vdots&\\

   (h_0,p_k) && 0 & 0 & 0 & -f_0 \quad -g_0 &0 \quad 0 &0 \quad 0 & f_0 \quad g_0 & 0 \quad 0 &\dots & 0 \quad 0 &
\end{array}
\]

Similar row-reductions can be done for the remaining hyperplanes $h_1$,...,$h_{|H|-1}$. After these row-reductions and permuting the rows so that the first $|H|$ rows are those that are non-zero in the columns corresponding to the hyperplanes, $M^p_S(P,H,I)$ can be reduced to the form
\[M^p_S(P,H,I)=
        \begin{pmatrix}
            {\rm Id}_{|H|} & *\\
            0 & B
        \end{pmatrix}
\]
where ${\rm Id}_H$ is the $|H| \times |H|$-identity matrix and $B$ is a $(2|P|-2) \times (2|P|-2)$ matrix. Note that $\det M^p_S(P,H,I)= \pm \det(B)$, so to find the pure condition it suffices to compute $\det(B)$. 

\begin{rmk}
    If the incidence geometry is a graph, and we are considering parallel redrawings in the plane, then $B$ is the \emph{parallel design matrix} of the graph - see \cite[Example 4.3.1]{discrete_matroids}.
\end{rmk}

The columns of $B$ come in pairs, with each pair indexed by an element of $\{p_0,\dots,p_{|P|-2}\}$.
Each row of $B$ is indexed by an incidence $(h_i,p_j)$ and all entries of such a row are $0$,
aside from the two entries corresponding to $p_j$, which will have values $f_i$ and $g_i$.
Since each point must lie on at least two hyperplanes,
this means that we can assume that every pair of columns corresponding to a single point $p_j$
will have a $2\times 2$ submatrix of the form
\[
    \begin{pmatrix}
        f_{i_1} & g_{i_1}\\
        f_{i_2} & g_{i_2}
    \end{pmatrix}.
\]
Let $l$ be the number of points that lie on exactly two hyperplanes
and without loss of generality assume that they are $p_0,\dots,p_{l-1}$.
Further assume without loss of generality that for $i = 0,\dots,l-1$,
the hyperplanes containing $p_i$ are $h_{2i}$ and $h_{2i+1}$.
Then, if the $i^{\rm th}$ row of $B$ has corresponding hyperplane $h_i$ and the $j^{\rm th}$ pair of columns has corresponding point $p_j$,
we have
\[
        B=
        \begin{pmatrix}
            B_0 & * & * & *\\
            0 & \ddots & * & * \\
            0 & 0 & B_{l-1} & * \\
            0 & 0 & 0 & B'
        \end{pmatrix}
    \]
\noindent
where each $B_i$ is a $2 \times 2$-block.
The $2 \times 2$-block $B_i$ in the columns corresponding to a point $p$ is of the form 
$B_i=
        \begin{pmatrix}
            \pm f_j & \pm g_j \\
            \pm f_k & \pm g_k \\
        \end{pmatrix}$
where $j$ and $k$ are such that $h_j$ and $h_k$ are incident to $p_i$.
Then $\det(B_i)=0$ means that the hyperplanes $h_j$ and $h_k$ are parallel
and that the bracket $[h_jh_k]=\det(B_i)$ will be a factor in $\det(M^S_p(P,H,I))$.
It remains to compute the determinant of the block $B'$ as a bracket polynomial.

\begin{ex}\label{Counterexample}

    \begin{figure}
    \[
        \begin{tikzpicture}
            \vertex (1) at (0,0)[label=left:$p_0$]{};
            \vertex (2) at (2,0)[label=right:$p_1$]{};
            \path (0,0) edge[thick, bend right] (2,0);
            \path (0,0) edge[thick, bend left] (2,0);
            \node at (1,0.75){$h_0$};
            \node at (1,-0.75){$h_1$};
        \end{tikzpicture}
    \]
    \caption{A basis of the $2$-plane matroid with no proper realizations.}
    \label{fig:2-plane_counterexample}
\end{figure}

    Let $S=(P,H,I)$ be the incidence geometry with two points $\{p_0,p_1\}$ and two lines $\{h_0, h_1\}$, that are both incident to both points $p_0$ and $p_1$ (Figure \ref{fig:2-plane_counterexample}. It is easy to verify that this incidence geometry is a basis of the $2$-plane matroid. Clearly, this incidence geometry does not have any
    realisations where the points and lines are distinct.
    
    The matrix $M_{p_0}^S(P,H,I)$ is

    \[
\begin{array}{ccccccccccccccccccccccccc}
     & & h_1 & h_2 & p_0 & p_1& & \\
   (h_{0},p_0) & \ldelim({6}{0.5em}& 1 & 0 & f_{0} \quad g_{0}& 0 \quad 0 & \rdelim){6}{0.5em} \\
   (h_{0},p_1) & & 1 & 0 & 0 \quad 0 & f_{0} \quad g_{0}\\
   (h_{1},p_0) & & 0 & 1 & f_{1} \quad g_{1}& 0 \quad 0 & \\
   (h_{1},p_1) & & 0 & 1 & 0 \quad 0 & f_{1} \quad g_{1}\\
    & & 0 & 0 & 1 \quad 0 & 0 \quad 0\\
    & & 0 & 0 & 0 \quad 1 & 0 \quad 0\\
\end{array}.
\]

Removing the row corresponding to the incidence $(h_0,p_0)$ from the row corresponding to $(h_0,p_1)$ and the row corresponding to $(h_1,p_0)$ from the row corresponding to $(h_1,p_1)$, gives the matrix

    \[
\begin{array}{ccccccccccccccccccccccccc}
     & & h_1 & h_2 & p_0 & p_1& & \\
   (h_{0},p_0) & \ldelim({6}{0.5em}& 1 & 0 & f_{0} \quad g_{0}& 0 \quad 0 & \rdelim){6}{0.5em} \\
   (h_{0},p_1) & & 0 & 0 & -f_0 \quad -g_0 & f_{0} \quad g_{0}\\
   (h_{1},p_0) & & 0 & 1 & f_{1} \quad g_{1}& 0 \quad 0 & \\
   (h_{1},p_1) & & 0 & 0 & -f_1 \quad -g_1 & f_{1} \quad g_{1}\\
    & & 0 & 0 & 1 \quad 0 & 0 \quad 0\\
    & & 0 & 0 & 0 \quad 1 & 0 \quad 0\\
\end{array}
\]
    which, after a permutation of the rows, is
    \[
\begin{array}{ccccccccccccccccccccccccc}
     & & h_1 & h_2 & p_0 & p_1& & \\
   (h_{0},p_0) & \ldelim({6}{0.5em}& 1 & 0 & f_{0} \quad g_{0}& 0 \quad 0 & \rdelim){6}{0.5em} \\
   (h_{1},p_0) & & 0 & 1 & f_{1} \quad g_{1}& 0 \quad 0 & \\
    & & 0 & 0 & 1 \quad 0 & 0 \quad 0\\
    & & 0 & 0 & 0 \quad 1 & 0 \quad 0\\
   (h_{0},p_1) & & 0 & 0 & -f_0 \quad -g_0 & f_{0} \quad g_{0}\\
   (h_{1},p_1) & & 0 & 0 & -f_1 \quad -g_1 & f_{1} \quad g_{1}\\
\end{array}
\]

    so the pure condition is simply the bracket $[h_0 h_1]$. 
\end{ex}

Example \ref{Counterexample} shows that not all bases of the $2$-plane matroid have proper realizations.
The $d$-uniform hypergraph with $d$ edges that all contain the same $d$ vertices are also bases of the $d$-plane matroid that have no proper realizations in dimension $d$.
One may now ask the following question: Does any basis of the $2$-plane matroid that does not contain the incidence geometry of Example \ref{Counterexample} as a subgeometry have proper realizations?

\begin{rmk}
    We will see in following examples that there are incidence geometries that do have proper realizations, but whose pure condition has factors $[h_i h_j]$, which is the pure condition of the incidence geometry with two points contained in two lines, $h_i$ and $h_j$. Having such a factor therefore does not distinguish which bases of the $2$-plane matroid have proper realisations.
\end{rmk}

\begin{ex} \label{Non-Fano}
    Consider the incidence geometry in Figure \ref{fig:Non-Fano}.
    We set set the pinned point $p$ to be $p_0$.
    Let $(f_i,g_i)$ denote the normal of the hyperplane $h_i$. Row-reducing gives two $2 \times 2$-blocks, $B_1=\begin{pmatrix}
     f_1 & g_1 \\
     f_5 & g_5
    \end{pmatrix}$
    and $B_2=\begin{pmatrix}
     f_2 & g_2 \\
     f_4 & g_4
    \end{pmatrix}$ on the diagonal. The pure condition will therefore have a factor $[h_1 h_5][h_2h_4]$, where $[h_i h_j]$ denotes the determinant of the matrix $\begin{pmatrix}
     f_i & g_i \\
     f_j & g_j
    \end{pmatrix}$. 
    The remaining block $B'$ has the following form:
\[B'=
        \begin{pmatrix}
        f_2 & g_2 &0 &0 & 0 &0 &0 & 0 \\
        f_3 & g_3 &-f_3 &-g_3 & 0 &0 &0 & 0 \\
        f_5 & g_5 &0 &0 & -f_5 & -g_5 &0 & 0 \\
        0 & 0 &f_0 &g_0& 0 &0 &0 & 0 \\
        0 & 0 &-f_3 &-g_3 & 0 &0 &f_3 & g_3 \\
        0 & 0 &0 &0 &-f_4 &-g_4 &f_4 & g_4 \\
        0 & 0 &0 &0 & f_0 &g_0 &0 &0 \\
        0 & 0 &0 &0 & 0 &0 &f_1 & g_1 \\
        \end{pmatrix}.
\]
Then $\det(B') = [h_2h_3] \det(A_1) + [h_2h_5] \det(A_2)$ where 
\[A_1=
        \begin{pmatrix}
        0 &0 & -f_5 & -g_5 &0 & 0 \\
        f_0 &g_0& 0 &0 &0 & 0 \\
        -f_3 &-g_3 & 0 &0 &f_3 & g_3 \\
        0 &0 &-f_4 &-g_4 &f_4 & g_4 \\
        0 &0 & f_0 &g_0 &0 &0 \\
        0 &0 & 0 &0 &f_1 & g_1 \\
        \end{pmatrix}
\]
and 
\[A_2=
        \begin{pmatrix}
        -f_3 &-g_3 & 0 &0 &0 & 0 \\
        f_0 &g_0& 0 &0 &0 & 0 \\
        -f_3 &-g_3 & 0 &0 &f_3 & g_3 \\
        0 &0 &-f_4 &-g_4 &f_4 & g_4 \\
        0 &0 & f_0 &g_0 &0 &0 \\
        0 &0 & 0 &0 &f_1 & g_1 \\
        \end{pmatrix}.
\]
$A_1$ can be further reduced to $\det(A_1) = [h_0 h_3] \det(A_1')$ where 
\[A_1'=
        \begin{pmatrix}
        -f_5 & -g_5 &0 & 0 \\
        f_0 &g_0 &0 &0 \\
        -f_4 &-g_4 &f_4 & g_4 \\
        0 &0 &f_1 & g_1 \\
        \end{pmatrix}.
\] Noting that $A_1'$ is block-lower triangular, we see that $\det(A_1)=[h_0h_3][h_0h_5][h_1h_4]$. Similarly, $\det(A_2) = [h_0h_4]\det(A_2')$, where \[A_2'=
        \begin{pmatrix}
        -f_3 &-g_3 &0 & 0 \\
        f_0 &g_0& 0 & 0 \\
        -f_3 &-g_3 & f_3 & g_3 \\
        0 &0 &f_1 & g_1 \\
        \end{pmatrix}.
\]
Noting that $A_2'$ is also block-lower triangular, we see that $\det(A_2)=-[h_0h_4][h_0h_3][h_1h_3]$. We see that $[h_0h_3]$ is a factor in both $A_1$ and $A_2$. 
In total, the pure condition of the incidence geometry in Figure \ref{fig:Non-Fano} becomes
\[
[h_1h_5][h_2h_4][h_0h_3] \cdot  ([h_2h_3][h_0h_5][h_1h_4]-[h_2h_5][h_0h_4][h_1h_3])=0.
\]
Notice that if any of the factors $[h_1h_5]$, $[h_2h_4]$ or $[h_0h_5]$ are 0, that means two of the hyperplanes must be parallel. As the relevant hyperplanes share a point, they must coincide. For the incidence geometry in Figure \ref{fig:Non-Fano}, two hyperplanes coinciding means that all hyperplanes coincide. Therefore, in any realizations of the incidence geometry where the hyperplanes do not coincide, the factor $[h_2h_3][h_0h_5][h_1h_4]-[h_2h_5][h_0h_4][h_1h_3]$ must be zero.
\end{ex}

\begin{figure}
\begin{center}
\begin{tikzpicture}

\draw[thick] (0,0) -- (4,0);
\draw[thick] (0,-2) -- (4,-3);
\draw[thick] (0,0)--(2,-2.5);
\draw[thick] (2,0) -- (4,-3);
\draw[thick] (4,0) -- (0,-2);

\draw[thick] (0,0) -- (4,-3);
\draw[thick] (2,0) -- (0,-2);
\draw[thick] (4,0) -- (2,-2.5);

\vertex (1) at (0,0)[label=above:$p_0$]{};
\vertex (2) at (2,0)[label=above:$p_1$]{};
\vertex (3) at (4,0)[label=above:$p_2$]{};

\vertex (4) at (0,-2)[label=below:$p_3$]{};
\vertex (5) at (2,-2.5)[label=below:$p_4$]{};
\vertex (6) at (4,-3)[label=below:$p_5$]{};

\vertex (4) at (0.88,-1.11)[label=left:$p_6$]{};
\vertex (5) at (1.6,-1.2)[label=above:$p_7$]{};
\vertex (6) at (2.91,-1.36)[label=above:$p_8$]{};

\node[anchor=south] at (1,0) {$h_0$};
\node[anchor=south] at (0,-0.8) {$h_1$};
\node[anchor=south] at (0.9,-0.8) {$h_2$};
\node[anchor=south] at (1.4,-0.8) {$h_3$};
\node[anchor=south] at (2.2,-0.8) {$h_4$};
\node[anchor=south] at (3,-0.8) {$h_5$};
\node[anchor=south] at (3.7,-0.8) {$h_6$};
\node[anchor=south] at (1,-2.7) {$h_7$};

\end{tikzpicture}
	\end{center}
	\caption{The non-Pappus configuration.}
	\label{Pappus_sub}
\end{figure}

\begin{ex}\label{ex:Pappus}
Consider the incidence geometry in Figure \ref{Pappus_sub}.
It consists of $9$ points and $8$ hyperplanes: $h_0=\{p_0,p_1,p_2\}$, $h_1=\{p_0,p_4,p_6\}$, $h_2=\{p_0,p_5,p_7 \}$ $h_3=\{p_1,p_3,p_6\}$, $h_4=\{p_1,p_5,p_8\}$, $h_5=\{p_2, p_3, p_7\}$, $h_6=\{ p_2, p_4,p_8\}$, and $h_7=\{p_3,p_4,p_5 \}$. Again, suppose that the point $p_0$ is pinned, and let $[h_i h_j]$ denote the determinant of the matrix  $\begin{pmatrix}
     f_i & g_i \\
     f_j & g_j
    \end{pmatrix}$ where $(f_i, g_i)$ and $(f_j, g_j)$ are the normals of the hyperplanes $h_i$ and $h_j$ respectively.

    Row-reducing the pinned matrix $M_{p_0}^S(P,H,I)$ gives three $2 \times 2$-blocks on the diagonal, namely $[h_1 h_3]$, $[h_2h_5]$ and $[h_4h_6]$, so the pure condition has a factor $[h_1 h_3][h_2 h_5][h_4h_6]$. The block $B'$ is the following matrix: 
    \[B'=
        \begin{pmatrix}
        f_1 & g_1 &0 &0 & 0 &0 &0 & 0 & 0 & 0 \\
        f_6 & g_6 &-f_6 &-g_6 & 0 &0 &0 & 0 & 0 & 0 \\
        f_7 & g_7 &0 &0 & 0 &0 &0 & 0 & -f_7 & -g_7 \\
        0 & 0 &-f_5 &-g_5 & 0 &0 &0 & 0 & f_5 & g_5 \\
        0 & 0 &f_0 &g_0 & 0 &0 &0 & 0 & 0 & 0 \\
        0 & 0 &0 &0 & 0 &0 &f_2 & g_2 & 0 & 0 \\
        0 & 0 &0 &0 & -f_3 &-g_3& 0 & 0 & f_3 & g_3 \\
        0 & 0 &0 &0 & -f_4 &-g_4 &f_4 & g_4 & 0 & 0 \\
        0 & 0 &0 &0 & 0 &0 &f_7 & g_7 & -f_7 & -g_7 \\
        0 & 0 &0 &0 & f_0 &g_0 &0 & 0 & 0 & 0 \\
        \end{pmatrix}.
\]

Similarly to in Example \ref{Non-Fano}, we can find that the determinant of $B'$ is
\[
[h_1h_7][h_0h_6]\cdot([h_0h_3][h_2h_4][h_5h_7]+[h_0h_4][h_2h_7][h_3h_5])+[h_1h_6][h_0h_5][h_0h_4][h_2h_7][h_3h_7].
\]
The pure condition of the incidence geometry in Figure \ref{Pappus_sub} is 
\[
    [h_1 h_3][h_2 h_5][h_4h_6]\cdot ([h_1h_7][h_0h_6]\cdot ([h_0h_3][h_2h_4][h_5h_7]+[h_0h_4][h_2h_7][h_3h_5])+[h_1h_6][h_0h_5][h_0h_4][h_2h_7][h_3h_7]). 
\] 
If any of the factors $[h_1 h_3]$, $[h_2 h_5]$ or $[h_4h_6]$ are zero, the hyperplanes $h_1$ and $h_3$, $h_2$ and $h_5$ or $h_4$ and $h_6$ must coincide.

Suppose that the bracket $[h_1 h_3]=0$, so that the hyperplanes $h_1$ and $h_3$ coincide. That forces the points $p_0$ and $p_1$ to coincide, as both points lie on the hyperplane $h_0$ and on the hyperplane $h_1=h_3$. Similarly, the points $p_3$ and $p_4$ have to coincide. That, in turn, forces the points $p_7$ and $p_8$ to coincide, as both points lie on the intersection of the hyperplane between $p_5$ and $p_0=p_1$ with the hyperplane between $p_2$ and $p_3=p_4$. See Figure \ref{degenerate_Pappus} for a realization where $[h_1 h_3]=0$. 

In any realization of the incidence geometry such that the points and hyperplanes are all distinct, it must hold that the following bracket polynomial vanishes
\[
    [h_1h_7][h_0h_6]\cdot([h_0h_3][h_2h_4][h_5h_7]+[h_0h_4][h_2h_7][h_3h_5])+[h_1h_6][h_0h_5][h_0h_4][h_2h_7][h_3h_7].
\]

\end{ex}
    \begin{figure}
	\begin{center}
	\begin{tikzpicture}
\filldraw[black] (2,0) circle (2pt);
\filldraw[black] (4,0) circle (2pt);

\filldraw[black] (2,-2.5) circle (2pt);
\filldraw[black] (4,-3) circle (2pt);
\filldraw[black] (2.91,-1.36) circle (2pt);

\draw[thick] (2,0) -- (4,0);
\draw[thick] (2,-2.5) -- (4,-3);

\draw[thick] (2,0) -- (4,-3);

\draw[thick] (4,0) -- (2,-2.5);

\draw[thick] (2,0)--(2,-2.5);
\filldraw[black] (2,-1.25) circle (2pt);

\node[anchor=south] at (2,0) {$p_0=p_1$};
\node[anchor=south] at (4,0) {$p_2$};

\node[anchor=north] at (2,-2.5) {$p_3=p_4$};
\node[anchor=north] at (4,-3) {$p_5$};

\node[anchor=south] at (3,0) {$h_0$};
\node[anchor=west] at (3.75,-0.5) {$h_5=h_6$};
\node[anchor=east] at (2,-0.5) {$h_1=h_3$};
\node[anchor=west] at (3.75,-2.5) {$h_2=h_4$};
\node at (3,-3) {$h_7$};

\node[anchor=west] at (2.91,-1.36) {$p_7=p_8$};
\node[anchor=east] at (2,-1.25) {$p_6$};
\end{tikzpicture}
	\end{center}
	\caption{A realization of the incidence geometry in Example \ref{ex:Pappus} with some coincident points and hyperplanes.}
	\label{degenerate_Pappus}
\end{figure}

\begin{figure}
    \[
        \begin{tikzpicture}
            \vertex (1) at (0,0)[label=left:$p_1$]{};
            \vertex (2) at (2,0)[label=right:$p_2$]{};
            \vertex (3) at (0,-2)[label=left:$p_3$]{};
            \vertex (4) at (2,-2)[label=right:$p_4$]{};
            \draw[thick, red] (0,0)--(2,0);
            \draw[thick, red] (0,0)--(0,-2);
            \draw[thick, red] (0,0)--(2,-2);
            \draw (2,0)--(2,-2);
            \draw (2,0)--(0,-2);
            \draw (2,-2)--(0,-2);
            \node at (1,0.5){$e_{12}$};
            \node at (-0.5,-1){$e_{13}$};
            \node at (2.5,-1){$e_{24}$};
            \node at (1,-2.5){$e_{34}$};
            \node at (0.75,-0.5){$e_{14}$};
            \node at (0.75,-1.5){$e_{23}$};
        \end{tikzpicture}
    \]
    \caption{$K_4$, with the tree $T$ in Example \ref{ex:treepoly} marked.}
    \label{fig:K_4-Tree}
\end{figure}

Let $G = (V,E)$ a graph that is a circuit of the $2d$-rigidity matroid, i.e.~$G$ satisfies $|E| = 2|V|-2$
and that if $V'$ is the set of vertices incident to some edge in $E' \subsetneq E$, then $|E'| \le 2|V| - 3$.
In this case, it follows from a result of Whiteley~\cite[Corollary 4.1.3]{discrete_matroids} that
if $(f_e,g_e)_{e \in E}$ is the vector of edge directions of a generic realization of $G$ in $\mathbb{R}^2$,
then the vector $(f_e/g_e)_{e \in E}$ lies in a hypersurface in $\mathbb{R}^E$ (by genericity, $g_e \neq 0$).
Martin~\cite{MartinGraphVarieties} calls this hypersurface the \emph{slope variety} of $G$
and gives a combinatorial method for finding a degree-$(|V|-1)$ polynomial whose vanishing defines this hypersurface;
he calls this polynomial the \emph{tree polynomial} of $G$.
We can view $G$ as an incidence geometry with point set $V$, hyperplane set $E$, with incidence set $I$ consisting of pairs $(e,v)$ whenever $e$ is incident to $v$.
Then, the slope variety of $G$ is also defined by taking the $2$-dimensional pure condition of $(V,E,I)$
and setting $g_e = 1$ for every $e \in E$.
This polynomial will also have degree $|V|-1$, and is thus equal to Martin's tree polynomial up to multiplication by a nonzero scalar.
We now illustrate this with an example.

\begin{ex}
\label{ex:treepoly}
    Let $S=K_4$ be the complete graph with vertex set $V=\{v_1,v_2, v_3, v_4\}$. Let $e_{ij}$ denote the edge $\{v_i,v_j\}$, let $E$ denote the set of all six such edges, and let $I=\{(e_{ij},v_k) \mid 1 \leq i < j \leq 4, k \in \{i,j\}\}$. It is easy to verify that $K_4=(V,E,I)$ is a basis of the $2$-plane matroid.
    The matrix $M^1_{K_4}(V,E,I)$ is shown at the top of Figure~\ref{fig: K_4 matrix}.
    Subtracting the row corresponding to the incidence $(e_{ij},v_i)$ from the row corresponding to the incidence $(e_{ij},v_j)$ for all $1\leq i < j \leq 4$ and permuting the rows gives the matrix shown in the middle of Figure~\ref{fig: K_4 matrix}.
    Hence the pure condition is the determinant of the $|E| \times (2|V|-2)$ matrix shown at the bottom of Figure~\ref{fig: K_4 matrix}.

    \begin{figure}
    \[
    \begin{array}{ccccccccccccccccccccccccc}
     & & e_{12} & e_{13} & e_{14} & e_{23} & e_{24} & e_{34} & v_1 & v_2 & v_3 & v_4& & \\
   (e_{12},v_1) & \ldelim({14}{0.5em} & 1 &0 &0 &0 &0 &0 &f_{12} \quad g_{14} & 0 \quad 0& 0 \quad 0 &0 \quad 0& \rdelim){14}{0.5em} \\
   (e_{12},v_2) & & 1 &0 &0 &0 &0 &0 &0 \quad 0 & f_{12} \quad g_{12}& 0 \quad 0 &0 \quad 0& \\
   (e_{13},v_1) & & 0 &1 &0 &0 &0 &0 &f_{13} \quad g_{13} & 0 \quad 0& 0 \quad 0 &0 \quad 0& \\
   (e_{13},v_3) & & 0 &1 &0 &0 &0 &0 &0 \quad 0 &0 \quad 0& f_{13} \quad g_{13} &0 \quad 0& \\
   (e_{14},v_1) & & 0 &0 &1 &0 &0 &0 &f_{14} \quad g_{14} & 0 \quad 0& 0 \quad 0 &0 \quad 0& \\
   (e_{14},v_4) & & 0 &0 &1 &0 &0 &0 &0 \quad 0 & 0 \quad 0& 0 \quad 0 &f_{14} \quad g_{14}& \\
   (e_{23},v_2) & & 0 &0 &0 &1 &0 &0 &0 \quad 0 & f_{23} \quad g_{23}& 0 \quad 0 &0 \quad 0& \\
   (e_{23},v_3) & & 0 &0 &0 &1 &0 &0 &0 \quad 0 & 0 \quad 0& f_{23} \quad g_{23} &0 \quad 0& \\
   (e_{24},v_2) & & 0 &0 &0 &0 &1 &0 &0 \quad 0 & f_{24} \quad g_{24}& 0 \quad 0 &0 \quad 0& \\
   (e_{24},v_4) & & 0 &0 &0 &0 &1 &0 &0 \quad 0 & 0 \quad 0& 0 \quad 0 &f_{24} \quad g_{24}& \\
   (e_{34},v_2) & & 0 &0 &0 &0 &0 &1 &0 \quad 0 & 0 \quad 0& f_{34} \quad g_{34} &0 \quad 0& \\
   (e_{34},v_4) & & 0 &0 &0 &0 &0 &1 &0 \quad 0 & 0 \quad 0& 0 \quad 0 &f_{34} \quad g_{34}& \\
   & & 0 &0 &0 &0 &0 &0 &1 \quad 0 & 0 \quad 0& 0 \quad 0 & 0 \quad 0& \\
     & & 0 &0 &0 &0 &0 &0 &0 \quad 1 & 0 \quad 0& 0 \quad 0 & 0 \quad 0& \\
    \end{array}
    \]
    \[
    \begin{array}{ccccccccccccccccccccccccc}
     & & e_{12} & e_{13} & e_{14} & e_{23} & e_{24} & e_{34} & v_1 & v_2 & v_3 & v_4& & \\
   (e_{12},v_1) & \ldelim({14}{0.5em} & 1 &0 &0 &0 &0 &0 &f_{12} \quad g_{12} & 0 \quad 0& 0 \quad 0 &0 \quad 0& \rdelim){14}{0.5em} \\
   (e_{13},v_1) & & 0& 1 &0 &0 &0 &0  &f_{13} \quad g_{13} & 0 \quad 0& 0 \quad 0 &0 \quad 0& \\
   (e_{14},v_1) & & 0 &0 &1 &0 &0 &0 &f_{14} \quad g_{14} & 0 \quad 0& 0 \quad 0 &0 \quad 0& \\
   (e_{23},v_2) & & 0 &0 &0 &1 &0 &0 &0 \quad 0 & f_{23} \quad g_{23}& 0 \quad 0 &0 \quad 0& \\
   (e_{24},v_2) & & 0 &0 &0 &0 &1 &0 &0 \quad 0 & f_{24} \quad g_{24}& 0 \quad 0 &0 \quad 0& \\
   (e_{34},v_3) & & 0 &0 &0 &0 &0 &1 &0 \quad 0 & 0 \quad 0& f_{34} \quad g_{34} &0 \quad 0& \\
   (e_{12},v_2) & & 0 &0 &0 &0 &0 &0 &-f_{12} \quad -g_{12} & f_{12} \quad g_{12}& 0 \quad 0 &0 \quad 0& \\
   (e_{13},v_3) & & 0 &0 &0 &0 &0 &0 &-f_{12} \quad -g_{12} &0 \quad 0& f_{13} \quad g_{13} &0 \quad 0& \\
   (e_{14},v_4) & & 0 &0 &0 &0 &0 &0 &-f_{14} \quad -g_{14} & 0 \quad 0& 0 \quad 0 &f_{14} \quad g_{14}& \\
   (e_{23},v_3) & & 0 &0 &0 &0 &0 &0 &0 \quad 0 & -f_{23} \quad -g_{23}& f_{23} \quad g_{23} &0 \quad 0& \\
   (e_{24},v_4) & & 0 &0 &0 &0 &0 &0 &0 \quad 0 & -f_{24} \quad -g_{24}& 0 \quad 0 &f_{24} \quad g_{24}& \\
   (e_{34},v_4) & & 0 &0 &0 &0 &0 &0 &0 \quad 0 & 0 \quad 0& -f_{34} \quad -g_{34} &f_{34} \quad g_{34}& \\
    & & 0 &0 &0 &0 &0 &0 &1 \quad 0 & 0 \quad 0& 0 \quad 0 & 0 \quad 0& \\
     & & 0 &0 &0 &0 &0 &0 &0 \quad 1 & 0 \quad 0& 0 \quad 0 & 0 \quad 0& \\
    \end{array}
\]
\[
    \begin{array}{ccccccccccccccccccccccccc}
     & & v_2 & v_3 & v_4& & \\
   e_{12} & \ldelim({6}{0.5em}& f_{12} \quad g_{12}& 0 \quad 0 &0 \quad 0& \rdelim){6}{0.5em} \\
   e_{13} & & 0 \quad 0& f_{13} \quad g_{13} &0 \quad 0& \\
   e_{14} & & 0 \quad 0& 0 \quad 0 &f_{14} \quad g_{14}& \\
   e_{23} & & -f_{23} \quad -g_{23}& f_{23} \quad g_{23} &0 \quad 0& \\
   e_{24} & & -f_{24} \quad -g_{24}& 0 \quad 0 &f_{24} \quad g_{24}& \\
   e_{34} & & 0 \quad 0& -f_{34} \quad -g_{34} &f_{34} \quad g_{34}& \\
\end{array}.
\]
    \caption{Matrices from Example~\ref{ex:treepoly}.
        The top matrix is $M^2_{K_4}(V,E,I)$.
        The middle matrix is obtained from the top matrix by row-reducing.
        Deleting the rows and columns of the middle matrix corresponding to the upper-left identity block yields
        the bottom matrix.
    }\label{fig: K_4 matrix}
\end{figure}

Let $[e_{ij}e_{kl}]$ denotes the determinant of the matrix
    $\begin{pmatrix}
     f_{ij} & g_{ij} \\
     f_{kl} & g_{kl}
    \end{pmatrix}$. 
It can be verified that the pure condition is 
\[[e_{12}e_{24}][e_{13}e_{23}][e_{14}e_{34}]-[e_{12}e_{23}][e_{13}e_{34}][e_{14}e_{24}]\] 

Dehomogenizing the polynomial by setting $g_{ij}=1$ for all $\{ij\} \in E$ gives that the pure condition is 
    \[(f_{12}-f_{24})(f_{13}-f_{23})(f_{14}-f_{34})-(f_{12}-f_{23})(f_{13}-f_{34})(f_{14}-f_{24})\]

which is the tree polynomial of $K_4$, as defined in \cite{MartinGraphVarieties}.
The dehomogenized pure condition of $K_4$ is the determinant of the matrix
\[
\begin{array}{ccccccccccccccccccccccccc}
     & & v_2 & v_3 & v_4& & \\
   e_{12} & \ldelim({6}{0.5em}& f_{12} \quad 1& 0 \quad 0 &0 \quad 0& \rdelim){6}{0.5em} \\
   e_{13} & & 0 \quad 0& f_{13} \quad 1 &0 \quad 0& \\
   e_{14} & & 0 \quad 0& 0 \quad 0 &f_{14} \quad 1& \\
   e_{23} & & -f_{23} \quad -1& f_{23} \quad 1 &0 \quad 0& \\
   e_{24} & & -f_{24} \quad -1& 0 \quad 0 &f_{24} \quad 1& \\
   e_{34} & & 0 \quad 0& -f_{34} \quad -1 &f_{34} \quad 1& \\
\end{array}.
\]

    Let $T$ be the spanning tree of $K_4$ with edges $e_{12}$, $e_{13}$ and $e_{14}$ (see Figure \ref{fig:K_4-Tree}). The unique circuit in $T \cup e_{23}$ is the cycle with edge-set $\{e_{12},e_{13},e_{23}\}$. Adding the row corresponding to $e_{12}$ to the row corresponding to $e_{23}$, and removing the row corresponding to $e_{13}$ from the row corresponding to $e_{23}$ gives the matrix

            \[
\begin{array}{ccccccccccccccccccccccccc}
   e_{12} & \ldelim({6}{0.5em}& f_{12} & 1& 0 & 0 &0 & 0& \rdelim){6}{0.5em} \\
   e_{13} & & 0 & 0& f_{13} & 1 &0 & 0& \\
   e_{14} & & 0 & 0& 0 & 0 &f_{14} & 1& \\
   e_{23} & & f_{12}-f_{23} & 0& f_{23}-f_{13} & 0 &0 & 0& \\
   e_{24} & & -f_{24} & -1& 0 & 0 &f_{24} & 1& \\
   e_{34} & & 0 & 0& -f_{34} & -1 &f_{34} & 1& \\
\end{array}.
\]

Using that the unique circuit in $T \cup e_{24}$ has edge-set $\{e_{12}, e_{14},e_{24}\}$ and that the unique circuit in $T \cup e_{34}$ has edge-set $\{e_{13}, e_{14},e_{34}\}$ and performing similar row-reductions, the matrix becomes

\[
\begin{array}{ccccccccccccccccccccccccc}
   e_{12} & \ldelim({6}{0.5em}& f_{12} & 1& 0 & 0 &0 & 0& \rdelim){6}{0.5em} \\
   e_{13} & & 0 & 0& f_{13} & 1 &0 & 0& \\
   e_{14} & & 0 & 0& 0 & 0 &f_{14} & 1& \\
   e_{23} & & f_{12}-f_{23} & 0& f_{23}-f_{13} & 0 &0 & 0& \\
   e_{24} & & f_{12}-f_{24} & 0& 0 & 0 &f_{24} -f_{14}& 0& \\
   e_{34} & & 0 & 0& f_{13}-f_{34} & 0 &f_{34} - f_{14} & 0& \\
\end{array}.
\]

Permuting columns now gives 

\[
\begin{array}{ccccccccccccccccccccccccc}
   e_{12} & \ldelim({6}{0.5em}& 1 & 0& 0 & f_{12} &0 & 0& \rdelim){6}{0.5em} \\
   e_{13} & & 0 & 1& 0 & 0 &f_{13} & 0& \\
   e_{14} & & 0 & 0& 1 & 0 & 0 & f_{14}& \\
   e_{23} & & 0 & 0& 0 & f_{12}-f_{23} & f_{23}-f_{13} & 0& \\
   e_{24} & & 0& 0& 0 & f_{12}-f_{24}  &0& f_{24} -f_{14}& \\
   e_{34} & & 0 & 0& 0 & 0 & f_{13}-f_{34}&f_{34} - f_{14} & \\
\end{array}
\]

so the dehomogenized pure condition is the determinant of the $(|V|-1) \times (|V|-1)$-matrix 

\[
\begin{array}{ccccccccccccccccccccccccc}
   e_{23} & \ldelim({3}{0.5em}&  f_{12}-f_{23} & f_{23}-f_{13} & 0&\rdelim){3}{0.5em} \\
   e_{24} & & f_{12}-f_{24}  &0& f_{24} -f_{14}& \\
   e_{34} & & 0 & f_{13}-f_{34}&f_{34} - f_{14} & \\
\end{array}
\]

which, by definition, is the tree polynomial with respect to $T$ (see Example 5.1 in \cite{MartinGraphVarieties}). 

In this example, pinning the vertex $v_1$ returned exactly the matrix that defines the tree polynomial. By Proposition \ref{prop: pure condition doesn't depend on choice of pinned point}, which vertex is pinned does not affect the determinant, up to multiplication by a non-zero scalar.
\end{ex}

\section{Overconstrained incidence geometries}
In the previous sections, we focus on incidence geometries that are bases of the $d$-plane matroid, since the pinned matrix $M_S^p(P,H,I)$ is square for such incidence geometries. However, there are incidence geometries that are not independent in the $d$-plane matroid, but that still have proper realizations in $d$-dimensional space.

Suppose that $(P,H,I)$ is an incidence geometry such that $|I| > |H|+d|P|-d$. Since $|I|>|H|+d|P|-d$, the rank of $M_S(P,H,I)$ is at most $d|P|+|H|-d$, assuming that the normals are generic. For this section, we will suppose that the rank of $M_S(P,H,I)$ for generic normals is $|H|+d|P|-d$. Pinning a point $p$ in the same way as in the last section then gives a matrix $M_S^p(P,H,I)$ of rank $|H|+d|P|$. The pinned matrix $M_S^p(P,H,I)$ will have at least one non-zero $(|H|+d|P|) \times (|H|+d|P|)$-minor. For any set of hyperplane normals such that $(P,H,I)$ has a nontrivial realization with those normals, $M_S^p(P,H,I)$ will have rank at most $|H|+d|P|-1$. For a set of feasible normals, the determinants of all $(|H|+d|P|) \times (|H|+d|P|)$-minors of $M_S^p(P,H,I)$ will need to be zero.

Considering the normals as variables, and setting the determinants of all $(|H|+d|P|) \times (|H|+d|P|)$-minors to zero gives a set of polynomial equations that need to be satisfied simultaneously. By Theorem \ref{rank_independent}, the rank of $M_S^p(P,H,I)$ does not depend on the pinned point, so the sets of normals for which the $(|H|+d|P|) \times (|H|+d|P|)$-minors vanish does not depend on the pinned point.
In practice however, there may be many non-zero $(|H|+d|P|) \times (|H|+d|P|)$-minors of $M_S^p(P,H,I)$, which makes computations impractical.
It would be interesting to find elegant descriptions of generating sets of the ideal generated by the $(|H|+d|P|) \times (|H|+d|P|)$ minors
of $M_S^p(P,H,I)$ in terms of the combinatorics of $(P,H,I)$.

\begin{figure}
\begin{center}
\begin{tikzpicture}

\draw[thick] (0,0) -- (4,0);
\draw[thick] (0,-2) -- (4,-3);
\draw[thick] (0,0)--(2,-2.5);
\draw[thick] (2,0) -- (4,-3);
\draw[thick] (4,0) -- (0,-2);

\draw[thick] (0,0) -- (4,-3);
\draw[thick] (2,0) -- (0,-2);
\draw[thick] (4,0) -- (2,-2.5);
\draw[thick] (0.88,-1.11)--(2.91,-1.36);

\vertex (1) at (0,0)[label=above:$p_0$]{};
\vertex (2) at (2,0)[label=above:$p_1$]{};
\vertex (3) at (4,0)[label=above:$p_2$]{};

\vertex (4) at (0,-2)[label=below:$p_3$]{};
\vertex (5) at (2,-2.5)[label=below:$p_4$]{};
\vertex (6) at (4,-3)[label=below:$p_5$]{};

\vertex (4) at (0.88,-1.11)[label=left:$p_6$]{};
\vertex (5) at (1.6,-1.2)[label=above:$p_7$]{};
\vertex (6) at (2.91,-1.36)[label=above:$p_8$]{};

\end{tikzpicture}
	\end{center}
	\caption{The Pappus configuration.}
	\label{Pappus}
\end{figure}

For example, let $(P,H,I)$ be the Pappus configuration, which is pictured in Figure \ref{Pappus}.
One can check using a symbolic computation package, such as e.g.~Sagemath~\cite{sage} or Macaulay2~\cite{M2},
that $M_S^p(P,H,I)$ has $324$ non-zero $(9+2 \times 9) \times (9+2 \times 9)$-minors.
Note that in the particular case of the Pappus configuration, the points $p_6$, $p_7$ and $p_8$ lie on a hyperplane in any realization of the configuration in Figure \ref{Pappus_sub} by Pappus theorem. So, to realize the Pappus configuration, it is both necessary and sufficient to find a realization of the configuration in Figure \ref{Pappus_sub}. 

\section{Acknowledgements}

The second author was supported by Knut and Alice Wallenberg Foundation, Sweden, Grant 2020.0007. 

\bibliographystyle{abbrv}
\bibliography{refs}

\end{document}